\documentclass{amsart}

\usepackage[foot]{amsaddr}

\author{Matthew D. Kvalheim}
\address
{University of Pennsylvania,
	Philadelphia, PA, USA}	
	
\email{kvalheim@seas.upenn.edu}

\makeatletter
\@namedef{subjclassname@2020}{%
  \textup{2020} Mathematics Subject Classification}
\makeatother   
\title[A generalization of the Hopf degree theorem]{A generalization of the Hopf degree theorem}
\usepackage{import} 

\usepackage[utf8]{inputenc}
\usepackage[T1]{fontenc}
\usepackage{lmodern}
\usepackage{amsmath}
\usepackage{amssymb}
\usepackage{amsthm}
\usepackage{stmaryrd} 
\usepackage{mathtools}
\usepackage{tikz-cd}
\usepackage{subfig}

\usepackage{thmtools}
\usepackage{thm-restate} 

\newcommand{\R}{\mathbb{R}}

\newcommand{\sph}{\mathbb{S}}

\theoremstyle{definition}

\newtheorem{Th}{Theorem}

\newtheorem*{Quest-non}{Question}

\newtheorem*{Th-non}{Theorem}

\newcommand{\thistheoremname}{}
\newtheorem*{genericthm}{\thistheoremname}
{\renewcommand{\thistheoremname}{Theorem~\ref{#1}$'$}%
	\begin{genericthm}}
	{\end{genericthm}}

\newtheorem*{Def*}{Definition}

\newtheorem{Rem}{Remark}

\usepackage{marvosym}
\usepackage[
hypertexnames,%
citecolor=blue,%
colorlinks=true,%
linkcolor=red%
]{hyperref}
\usepackage[all]{hypcap}

\begin{document}
	
	\begin{abstract}	
    The Hopf theorem states that homotopy classes of continuous maps from a closed connected oriented smooth $n$-manifold $M$ to the $n$-sphere are classified by their degree.
    Such a map is equivalent to a section of the trivial $n$-sphere bundle over $M$.
    A generalization of the Hopf theorem is obtained for sections of nontrivial oriented $n$-sphere bundles over $M$.
	\end{abstract}

\maketitle
	
   Throughout, let $M$ be a compact connected oriented smooth $n$-dimensional manifold without boundary.
   The Hopf degree theorem is as follows \cite[p.~51]{milnor1965topology}.
   
   \begin{Th-non}[Hopf]
   A pair of continuous maps $f,g\colon M\to \sph^n$ to the $n$-sphere are homotopic if and only if $\deg(f) = \deg(g)$.
   \end{Th-non}
   Homotoping $f$ and $g$ to smooth maps, denoted by the same symbols, does not affect their degrees.
   Let $c\in \sph^n$ be a regular value for both.
   Then $\deg(f) = I(f,c)$ and $\deg(g)=I(g,c)$ coincide with oriented intersection numbers.
   If instead $c$ is viewed as the constant map $M\to \{c\}\subset \sph^n$, rather than as a point in $\sph^n$, these can be viewed as oriented intersection numbers of maps $M\to \sph^n$ \cite{milnor1965topology,guillemin1974differential,hirsch1976differential}. 
   Such a map is equivalent to a section of the trivial sphere bundle $M\times \sph^n\to M$.
   This motivates the following question. 
\begin{Quest-non} 
What can be said when the sphere bundle is nontrivial and $M\to \{c\}$ is replaced by a general continuous section?
\end{Quest-non}   
A partial answer is given by the following generalization of the Hopf theorem, which is the special case that $\pi$ is a trivial bundle and $Z$ is constant.
   \begin{Th}\label{th:gen-Hopf}
   Let $X,Y,Z$ be continuous sections of a smooth $n$-sphere bundle $\pi\colon E\to M$, where $E$ is oriented.
   Then $X$, $Y$ are homotopic through sections if and only if $I(X,Z)=I(Y,Z)$.
   \end{Th}  
   
   \begin{Rem}
   The proof works by reducing to the situation of the Hopf theorem.
   \end{Rem}
   \begin{proof}
   Homotopy invariance of the oriented intersection number implies that $I(X,Z)=I(Y,Z)$ if $X$ and $Y$ are homotopic.
   
   Conversely, assume that $I(X,Z)=I(Y,Z)$.
   It suffices to show that $X$ and $Y$ are homotopic through sections to some common other section.
   By approximation techniques we may assume (after preliminary homotopies through sections) that $X$, $Y$, $Z$ are smooth and transverse (cf. \cite[p.~56, Ex.~3]{hirsch1976differential}).
   In particular, $X$, $Y$  intersect $Z$ only at finitely many points in $M$.
   
   Finiteness and connectedness of $M$ imply that these intersection points are contained in the interior of a compact set $B\subset M$ diffeomorphic to a ball in $\R^n$ \cite{michor1994n}.
   Since $B$ is smoothly contractible there is a fiber-preserving diffeomorphism
   $$\pi^{-1}(B)\approx B \times \sph^n$$
   with respect to which sections over $B$ may be viewed as $\sph^n$-valued and $Z|_B$ is constant.
   After homotopies through sections we may assume that 
   \begin{equation}\label{eq:coincide}
   \textnormal{$X|_B$ and $Y|_B$ coincide with $-Z|_B$ on a neighborhood of $\partial B$ in $B$.}
   \end{equation}
   
   It follows that $X$, $Y$ induce, by collapsing the boundary of $B$, self-maps
   $$\sph^n\approx B/\partial B \to \sph^n$$
   whose degrees coincide with $I(X,Z)=I(Y,Z)$.
   The Hopf theorem and its proof imply that the $X$-induced self-map is homotopic to the $Y$-induced self-map through maps sending $[\partial B]\in \sph^n$ to the constant $-Z|_B\in \sph^n$ \cite[pp.~50--51]{milnor1965topology}. 
   By \eqref{eq:coincide} this homotopy extends by the constant one to yield a global homotopy of $X$, $Y$ through sections.
   \end{proof}
   
   \begin{Rem}[cf. {\cite[p.~50, Remarks]{milnor1965topology}}]
   Theorem~\ref{th:gen-Hopf} easily generalizes to let $M$ have a nonempty boundary $\partial M$: 
   continuous sections $X$, $Y$ disjoint from $Z|_{\partial M}$ are homotopic through sections disjoint from $Z|_{\partial M}$  if and only if $I(X,Z)=I(Y,Z)$.
   \end{Rem}

   	\bibliographystyle{amsalpha}
   	\bibliography{ref}

\providecommand{\bysame}{\leavevmode\hbox to3em{\hrulefill}\thinspace}
\providecommand{\MR}{\relax\ifhmode\unskip\space\fi MR }
\providecommand{\MRhref}[2]{%
  \href{http://www.ams.org/mathscinet-getitem?mr=#1}{#2}
}
\providecommand{\href}[2]{#2}
\begin{thebibliography}{MV94}

\bibitem[GP10]{guillemin1974differential}
V~Guillemin and A~Pollack, \emph{Differential topology}, AMS Chelsea
  Publishing, Providence, RI, 2010, Reprint of the 1974 original. \MR{2680546}

\bibitem[Hir94]{hirsch1976differential}
M~W Hirsch, \emph{Differential topology}, Graduate Texts in Mathematics,
  vol.~33, Springer-Verlag, New York, 1994, Corrected reprint of the 1976
  original. \MR{1336822}

\bibitem[Mil97]{milnor1965topology}
J~W Milnor, \emph{Topology from the differentiable viewpoint}, Princeton
  Landmarks in Mathematics, Princeton University Press, Princeton, NJ, 1997,
  Based on notes by David W. Weaver, Revised reprint of the 1965 original.
  \MR{1487640}

\bibitem[MV94]{michor1994n}
P~W Michor and C~Vizman, \emph{n-transitivity of certain diffeomorphism
  groups}, Acta Math. Univ. Comenianae \textbf{63} (1994), no.~2, 221--225.

\end{thebibliography}

   \section*{Acknowledgments}
   This work was funded by ONR N00014-16-1-2817, a Vannevar Bush Faculty Fellowship  sponsored by the Basic Research Office
   of the Assistant Secretary of Defense for Research and Engineering.

\end{document}